\documentclass[
11pt,%
tightenlines,%
twoside,%
onecolumn,%
nofloats,%
nobibnotes,%
nofootinbib,%
superscriptaddress,%
noshowpacs,%
centertags]%
{article}
\usepackage{amsmath,amssymb}

\newtheorem{definition}{Definition}
\newtheorem{theorem}{Theorem}
\newtheorem{corollary}{Corollary}

\begin{document}

\title{Degrees of selector functions and relative computable categoricity}
\author{I. Sh. Kalimullin\footnote{The work was supported by the Theoretical Physics and Mathematics
Advancement Foundation ``BASIS''.}}

\maketitle
\abstract{\noindent We study the degrees of selector functions related to the degrees in which a rigid computable structure is relatively computably categorical. It is proved that for some structures such degrees can be represented as the unions of upper cones of c.e. degrees. In addition we show that there are non-c.e. upper cones realized as the degrees in which some computable structure is relatively computably categorical.}

\

\noindent {\bf Introduction.} We start with the main definition of the notion studied through the paper.
\begin{definition} Let $\{C_i\}_{i\in\omega}$  be a  computable sequence of  differences of c.e. sets, $C_i=A_i\setminus B_i$. A  {\em selector} function for $\{C_i\}_{i\in\omega}$ is any function $f$ such that $f(i)\in C_i$ for every $i\in\omega$.
\end{definition}

As an example for a computable  rigid structure $\mathcal A$ (on the domain $\omega$) consider
 the sequence $\{C_i=A_i\setminus B_i\}_{i\in\omega}$, where
$A_i$ and $B_i$ are the following c.e. sets of existential formulae in the language of $\mathcal A$:
$$
A_i=\{\Phi\mid \mathcal A\models \Phi(i)\},
$$
$$
B_i=\{\Phi\mid \mathcal A\models (\exists j\ne i) \Phi(j)\}.
$$
Then an existence of $\mathbf x$-computable selector function for $\{C_i\}_{i\in\omega}$ implies the relative $\mathbf x$-computable categoricity of the structure $\mathcal A$, i.e., for every isomporphic copy  $\mathcal B\cong \mathcal A$ there is a $(\deg(\mathcal B)\cup\mathbf x)$-computable isomorphism from $\mathcal B$ onto $\mathcal A$. Moreover,
by the well-known result of 
 Ash, Knight, Manasse and Slaman~\cite{AKMS} and Chisholm~\cite{Ch}
the inverse also holds:  if  a rigid computable structure $\mathcal A$ is $\mathbf x$-computably categorical
for some degree $\mathbf x $ then for some finite constant enrichment $\widetilde{\mathcal A}=(\mathcal A,\vec a)$ the sequence $\{\widetilde C_i=\widetilde A_i\setminus\widetilde B_i\}_{i\in\omega}$ defined as above
will have an $\mathbf x$-computable selector function.

It is easy to check here that the degrees of such selector functions does not actually depend on the choice of the constants $\vec a$. Moreover, due a possible conjuctions of existential formulae the degrees of the selector functions in this case are the same as the degrees of the weak selector functions defined below.

\begin{definition} Let
$\{D_n\}_{n\in\omega}$ be the  standard (canonical) numbering of all finite subsets of $\omega$ defined via $n=\sum_{x\in D_n}2^x.$
  A  {\em weak selector} function for a  computable sequence  of  differences of c.e. sets  $\{C_i=A_i\setminus B_i\}_{i\in\omega}$
 is a selector function for the sequence $\{\widehat C_i=\widehat A_i\setminus \widehat B_i\}_{i\in\omega}$, where
$$
\widehat A_i=\{n\mid D_n\subseteq A_i\}\mbox { and }
\widehat B_i=\{n\mid D_n\subseteq B_i\}.
$$
\end{definition}

\begin{theorem}\label{iff}
\begin{enumerate} \item
For every rigid computable structure  $\mathcal A$ there is a   computable sequence  of  differences of c.e. sets
$\{C_i\}_{i\in\omega}$ such that for each degree $\mathbf x$ the structure
$\mathcal A$ is relatively $\mathbf x$-computably categorical if and only if there is an $\mathbf x$-computable weak selector function for $\{C_i\}_{i\in\omega}$.
\item
For every computable sequence  of  differences of c.e. sets
$\{C_i\}_{i\in\omega}$ there is
a rigid computable structure  $\mathcal A$
 such that for each degree $\mathbf x$ the structure
$\mathcal A$ is relatively $\mathbf x$-computably categorical if and only if there is an $\mathbf x$-computable weak selector function for $\{C_i\}_{i\in\omega}$.
\end{enumerate}
\end{theorem}

The next theorem shows that $1$-generic oracles can not compute non-trivial selector functions.
\begin{theorem}
If a degree $\mathbf x$ is $1$-generic and there is an $\mathbf x$-computable (weak) selector function for a computable sequence  of  differences of c.e. sets
then there is
a computable (weak) selector function for this sequence.
\end{theorem}
\begin{corollary}
If a computable rigid structure $\mathcal A$  is relatively $\mathbf x$-computably categorical for a $1$-generic degree $\mathbf x$ then $\mathcal A$ is  relatively computably categorical.
\end{corollary}

In contrast with Theorem 2 an interesting example of computable sequence  of  differences of c.e. sets $\{C_i\}_{i\in\omega}$  with nontrivial properties of selector functions appear if we define  for any pair of c.e. sets $U\subseteq V$
$$C_i=\begin{cases}
\{0\}, & \mbox{if }i\notin V;\cr
\{0,1\}, & \mbox{if }i\in V\setminus U;\cr
\{1\}, & \mbox{if }i\in U.
\end{cases}$$
This can be defined also via $C_i=A_i\setminus B_i$, where $$A_i=\{0\}\cup\{1\mid i\in V\},$$  $$B_i=\{0\mid i\in U\}.$$
Then it is easy to see that an $\mathbf x$-computable  selector function for  $\{C_i\}_{i\in\omega}$ exists if and only if there is an $\mathbf x$-computable set $X$ such that  $U\subseteq X\subseteq  V.$ It is clear also that
$$\{0,1\}\subseteq A_i\implies 1\in C_i=A_i\setminus B_i,$$
so that
an existence of $\mathbf x$-computable selector function for such  $\{C_i\}_{i\in\omega}$  is equivalent to an existence of
$\mathbf x$-computable weak selector function for $\{C_i\}_{i\in\omega}$.

In the particular case $U=V$ we easily can build a computable sequence  of  differences of c.e. sets with a unique selector function of the same c.e. degree as $U$ and $V$. We can extend this as the following.
\begin{theorem}\label{union}
For every finite sequence of c.e. degrees $\mathbf a_1, \mathbf a_2,\dots, \mathbf a_k$ there is a pair of c.e. sets $U\subseteq V$ such that an $\mathbf x$-computable set $X$,  $U\subseteq X\subseteq  V,$ exists if and only if
$$
 \mathbf a_1\le \mathbf x\mbox{ or } \mathbf a_2\le \mathbf x\mbox{ or }\dots \mbox{ or } \mathbf a_k\le \mathbf x.
$$
\end{theorem}
By Theorem~\ref{iff} we can code the degrees of  sets $X$ such that  $U\subseteq X\subseteq  V$ into the degrees in which a computable structure is relatively computably categorical.
\begin{corollary}
For every finite sequence of c.e. degrees $\mathbf a_1, \mathbf a_2,\dots, \mathbf a_k$ there is a computable rigid structure $\mathcal A$ such that
$\mathcal A$ is relatively $\mathbf x$-computably categorical if and only if
$$
 \mathbf a_1\le \mathbf x\mbox{ or } \mathbf a_2\le \mathbf x\mbox{ or }\dots \mbox{ or } \mathbf a_k\le \mathbf x.
$$
\end{corollary}

Theorem~\ref{union} allows to build more non-trivial and more non-uniform examples of computable sequence  of  differences of c.e. sets, e.g., we can apply Theorem~\ref{union} for $\mathbf a_0>\mathbf a_1.$
By this way we can make only c.e. degrees as the least degrees of selector functions. The following theorem allows to find the least degrees of selector functions among $2$-CEA non-c.e. degrees.


\begin{theorem} Let $\mathbf e$ be a c.e. degree, and let $F\in\mathbf f$ be an $\mathbf e$-c.e. set such that there is a $\Delta^0_2$-approximation $F(x)=\lim_sF_s(x)$, $F_s(x)\in\{0,1\}$, with the property
$$F_s(x)\ne F_{s+1}(x) \implies F_s(y)=F_t(y)$$
for all $x<y<s<t$. Then there is a pair of c.e. sets $U\subseteq V$ such that an $\mathbf x$-computable set $X$,  $U\subseteq X\subseteq  V,$ exists if and only if   $\mathbf e\cup\mathbf f\le \mathbf x.$
\end{theorem}

The degrees $\mathbf e\cup\mathbf f$ satisfying the conditions of Theorem 4 form sufficiently large class of $2$-CEA degrees. Indeed, Jockusch and Shore~\cite{JS} have 
constructed a $\Delta^0_2$ $2$-CEA degree not belonging to any given uniform  $\Delta^0_2$ class. We can note that the construction of the corresponding set $F\in\mathbf f$ assumes only one witness per requirement, and the requirements are satisfied in a finite injury priority manner. If  a witness $x$ enters or leaves $F$ during the construction at a stage $s$ then other assigned earlier witnesses $y>x$ of lower priority are initialized at this stage, and hence  $F_s(y)= F_t(y)$ for $t>s.$. Therefore, the condition
$$x<y<s<t\;\&\; F_s(x)\ne  F_{s+1}(x) \implies  F_s(y)= F_t(y)$$
from Theorem 4  holds, so that we have proven the following statement. 

\begin{corollary}\label{2cea}
For every uniform  $\Delta^0_2$ class of degrees $\mathcal C$ (e.g., $\mathcal C=$ the c.e. degrees, $\mathcal C=$ the $2$-c.e. degrees, etc.) there are a $2$-CEA degree $\mathbf f\notin \mathcal C$ and a computable rigid structure $\mathcal A$ such that
$\mathcal A$ is relatively $\mathbf x$-computably categorical if and only if $\mathbf f\le \mathbf x$.
\end{corollary}

If $\mathcal C=$ the c.e. degrees the proof produces a non-c.e. $2$-c.e. degree $\mathbf f$ 
If $\mathcal C=$ the $2$-c.e. degrees by
the result of  Arslanov, LaForte and Slaman~\cite{ALS} we can not  produce a $3$-c.e. degree $\mathbf f\notin \mathcal C$. But the proof of Theorem~4 allows to repeatedly apply the arguments through the $n$-CEA hierararchy. Namely the proof of Theorem~4 can be adapted for $\mathbf e=\mathbf g_i$ and $\mathbf f=\mathbf g_{i+1}$, if
 $$\mathbf g_1<\mathbf g_2<\mathbf g_3<\cdots$$
and 
 $$\mathbf g_1\mbox { is a c.e. degree}, \mathbf g_2\mbox { is a $\mathbf g_1$-c.e. degree}, \mathbf g_3\mbox { is a $\mathbf g_2$-c.e. degree},\dots$$
 Here we need only check that each $\mathbf g_i$ contains an appropriate set $G_i\in \mathbf g_i$ with a $\Delta^0_2$-approximation satisfying the condition from Theorem~4. The following theorem can be proved by this way but we prefer to give a direct proof, where an appropriate approximation is needed only for the final degree $\mathbf f=\mathbf g_n.$
 
\begin{theorem} For every $n>1$ there are an $n$-c.e. degree $\mathbf f$ which is not $(n-1)$-c.e. and  a pair of c.e. sets $U\subseteq V$ such that an $\mathbf x$-computable set $X$,  $U\subseteq X\subseteq  V,$ exists if and only if   $\mathbf f\le \mathbf x.$
\end{theorem}
\begin{corollary}\label{2cea}
For every $n>1$ there are a non-$(n-1)$-c.e. $n$-c.e. degree $\mathbf f$ and  a computable rigid structure $\mathcal A$ such that
$\mathcal A$ is relatively $\mathbf x$-computably categorical if and only if $\mathbf f\le \mathbf x$.
\end{corollary}

The rest of the paper is devoted to the proofs of the theorems above. We use the monograph~\cite{So} as the source of  used notations and terminology.
 
\section{The proof of Theorem 1}

{\em The proof of Part 1.} By the result of  Ash, Knight, Manasse and Slaman~\cite{AKMS} and Chisholm~\cite{Ch} a rigid computable structure is $\mathcal A$ is relatively $\mathbf x$-computably categorical if and only if for some tuple $\vec a$ from $\mathcal A$ there is a $\mathbf x$-computable enumeration of existential formulae  defining all individual elements in the $(\mathcal A,\vec a)$.

If this never happen then we can simply define the sequence  $C_i=\emptyset$ which obviously has no selector function.

Otherwise, if for some $\vec a$ there is a collection of  existential formulae defining all individual elements in the $(\mathcal A,\vec a)$, then the  problem of enumeration such a collection does not depend on the choice of $\vec a$. Therefore, for a fixed $\vec a$ we can define  the sequence $\{C_i=A_i\setminus B_i\}_{i\in\omega}$ with
$$
A_i=\{\Phi\mid (\mathcal A,\vec a)\models \Phi(i)\},
$$
$$
B_i=\{\Phi\mid (\mathcal A, \vec a)\models (\exists j\ne i) \Phi(j)\}.
$$
%

\

\noindent {\em The proof of Part 2.} Suppose a computable sequence  of  differences of c.e. sets $\{C_i\}_{i\in\omega}$, $C_i=A_i\setminus B_i\ne \emptyset,$ be given. Without loss of generality we can assume that $B_i\ne\emptyset,$ $B_i\subseteq A_i$, and $A_i\cap A_j=\emptyset$ for $i\ne j.$ If not, we can consder the sequence  $\{\widetilde C_i\}_{i\in\omega}$, $\widetilde C_i=\widetilde  A_i\setminus \widetilde  B_i$ with the same degrees of selector functions:
$$
\widetilde A_i=\{\langle i,0\rangle\}\cup \{\langle i,x+1\rangle\mid x\in A_i\},
$$
$$
\widetilde B_i=\{\langle i,0\rangle\}\cup \{\langle i,x+1\rangle\mid x\in A_i\cap B_i\}.
$$
Let $a$ and $b$ be injective computable functions such that ${\rm rng}\; a=\cup_{i\in\omega}A_i$ and ${\rm rng}\; b=\cup_{i\in\omega}B_i$. Also let $h$ be a computable function such that $h(i)\in B_i$ for every $i\in \omega.$

We construct a structure $\mathcal A$ on the domain $\omega$ in the in the language of infinitely many unary functions $e_0,e_1,e_2,\dots,$ by the following:
$$
e_k(n)=\begin{cases} 4i,&\mbox{if }n=4j+2\mbox{ and }a(j)=k\in A_i;\cr
 4i+1,&\mbox{if }n=4j+3\mbox{ and }b(j)=k\in B_i;\cr
n,&\mbox{otherwise.}
\end{cases}
$$
Despite the infinity of the functional language the structure $\mathcal A$ is locally finite so that we can find a computable function $c$ such that the finite set $D_{c(s)}$ is the set generated from the elements $0,1,\dots,s$.

Suppose that  there is there is an $\mathbf x$-computable weak selector function $f$ for $\{C_i\}_{i\in\omega}$, i.e.,
$$
D_{f(i)}\subseteq A_i\mbox{ and } D_{f(i)}\not\subseteq B_i \mbox { for every }i\in\omega.
$$
Then we can build an $\mathbf x$-computable list of existential formulae defining all elements of $\mathcal A:$
$$
\Phi_{4i}(x)=\mbox{   }\mathop{\mbox{\Large \&}}_{k\in D_{f(i)}}(\exists z\ne x) [e_k(z)=x];
$$
$$
\Phi_{4i+1}(x)= \mbox{   }(\exists y\ne x)\Phi_{4i}(y)\;\&\;(\exists z\ne x )[e_{h(i)}(z)=x];
$$
$$
\Phi_{4i+2}(x)= \mbox{   }\Phi_{4i}(e_{a(i)}(x));
$$
$$
\Phi_{4i+3}(x)= \mbox{   }\Phi_{4i+1}(e_{b(i)}(x)).
$$
By~\cite{AKMS} the structure $\mathcal A$ is relatively $\mathbf x$-computably categorical.

Conversely, suppose the structure $\mathcal A$ is relatively $\mathbf x$-computably categorical. Then by ~\cite{AKMS} for some tuple $\vec a$ from $\mathcal A$ there is an $\mathbf x$-computable sequence $\{\Phi_i(x)\}_{i\in\omega}$ of existential formulae such that
$$i=j\iff (\mathcal A,\vec a)\models \Phi_i(j).$$
Fix an $s_*\in\omega$ such that  $\vec a\subseteq D_{c(s_*)}$. Then for every $i$ we can $\mathbf x$-computably find an $s_i>s_*$ such that  $4i\in D_{c(s_i)}$ and
$$(D_{c(s_i)},e_0, e_1, e_2,\dots, \vec a)\models \Phi_{4i}(4i).$$
Now if $4i\notin D_{c(s_*)}$ then an $\mathbf x$-computable weak selector function $f$ for  $\{C_i\}_{i\in\omega}$ can be defined for  as follows:
$$
D_{f(i)}=\{k\mid (\exists n\in D_{c(s_i)})[n\ne 4i\;\&\; e_x(n)=4i]\}.
$$
It is easy to see that $D_{f(i)}\subseteq A_i$ and $D_{f(i)}\not\subseteq B_i$ since otherwise we would have $ (\mathcal A,\vec a)\models \Phi_{4i}(4i+1)$. It is enough now to apppropriately extend the definition of $f$ for finitely many $i$ with  $4i\in D_{c(s_*)}$.

 \section{The proof of Theorem 2}

If $\mathbf x$ is $1$-generic then there is a set  $X\in\mathbf x$ such that for every c.e. set $W\subseteq 2^{<\omega}$ there is a string $\sigma\subset X$ with the property
$$\sigma\in W\mbox{ or }(\forall\tau\supseteq\sigma)[\tau\notin W].$$

Suppose that $f=\{e\}^X$ is an $\mathbf x$-computable selector function   for a  computable sequence  of  differences of c.e. sets  $\{C_i=A_i\setminus B_i\}_{i\in\omega}$. Let
$$W=\{\sigma\in 2^{<\omega}\mid (\exists i)[\{e\}^\sigma(i)\downarrow\in B_i]\}.$$
Since $f(i)=\{e\}^X(i)\in A_i\setminus B_i$ we can not have $\sigma\subset X$. Since $W$ is c.e. there is an $\sigma\subset X$ such that $\tau\notin W$  for all $\tau\supseteq\sigma$.

Then we can find a computable selector function  for  $\{C_i=A_i\setminus B_i\}_{i\in\omega}$: let $g(i)=\{e\}^{\sigma_i}(i)$, where $\sigma_i\supseteq\sigma$ is the first found string such that $\{e\}^{\sigma_i}(i)\in A_i$.

The statement of Theorem 2 for weak selector functions also follows from the arguments above since the weak selector functions for $\{C_i=A_i\setminus B_i\}_{i\in\omega}$ are the selector function for the modified sequence   $\{\widehat C_i=\widehat A_i\setminus \widehat B_i\}_{i\in\omega}$,
$$
\widehat A_i=\{n\mid D_n\subseteq A_i\}\mbox { and }
\widehat B_i=\{n\mid D_n\subseteq B_i\}.
$$

\section{The proof of Theorem 3} Without loss of generality we can assume that the given c.e. degrees  $\mathbf a_1, $ $\mathbf a_2,\dots,$ $ \mathbf a_k$ are all non-zero.
For $1\le i\le k$ let $A_i$ be a c.e. set such that $A_i\in \mathbf a_i$. 

To proceed we need the following trick similar to the Dekker's deficiency  set: if  $A={\rm rng\;} a$  and $V={\rm rng\;} v$ are infinite c.e. sets,  where $a$ and $v$ are computable injective functions, then we can consider the c.e. subset of $V:$
$$
A^V=\{v(s)\mid (\exists t>s )[a(t)<v(s)]\}\subseteq V.
$$
It is easy to see that for every $X\subseteq\omega$
$$A^V\subseteq^* X\implies V\subseteq^*X\mbox{ or } A\le_T X,$$
since every element  $v(s)\in V\setminus (A^V\cup X)=^* V\setminus X$ gives a computation of $A(x)$ for each $x<v(s)$:
$$x\in A\iff x\in\{a(0),a(1),\dots, a(s)\}.$$
In particular, if $A$ is not computable then the c.e. set $A^V$ is infinite. Also we have $A^V\le_T A$ since $A\restriction x\subseteq\{a(0),a(1),\dots, a(n)\}$ implies
$$
x\in A^V\iff(\exists t\le n) (\exists s<t)[x=v(s)>a(t)].
$$

Since each c.e. set $A_i\in\mathbf a_i$ is not computable  we can now consider the chain of infinite c.e. sets
$$U=V_k\subseteq V_{k-1}\subseteq\cdots\subseteq V_2\subseteq V_1=V,$$
where $V_1=A_1,$ and
$V_{i+1}=A_{i+1}^{V_i}$ for $i<k.$

By the arguments from above we have $V_i\le_T A_i$ for each $i$.  Also, if  $U\subseteq X\subseteq V$ for a set $X$ then  we have $A_i\le_T X$ where $i\le k$ is the least index such that $V_i\subseteq^* X$.

\section{The proof of Theorem 4}

 For a fixed c.e. set $E\in\mathbf e$ we fix an index $e$ such that $F=W_e^E={\rm dom\;}\{e\}^E$.
We re-define the
 $\Delta^0_2$-approximation $F(x)=\lim_sF_s(x)$  for $s=0,1$ by setting $F_0(x)=1$ and $F_1(x)=0$ for each $x$.
%

Note that the property
$$x<y<s<t\;\&\; F_s(x)\ne  F_{s+1}(x) \implies  F_s(y)= F_t(y)$$
again holds for the modified approximation (just because there are no tuples $x<y<s\le 1$).
Then for the set
$$
\widetilde F=\{\langle y,t\rangle\mid  F_{t}(y)= F(y)=1\;\&\; F_{t+1}(y)=0\}
$$
we have $ F\equiv_T \widetilde F$ due $y\in F\iff \langle y,0\rangle\in\widetilde F.$ Also we have  $\widetilde U\subseteq \widetilde F\subseteq \widetilde V$ for the c.e. sets

\medskip

$
\widetilde U=\{\langle y,t\rangle\mid
(\exists x)(\exists s)[x{<}y{<}s\;\&\;  F_s(x)\ne  F_{s+1}(x)\;\&\;  F_s(y)= F_{t}(y)=1]
$

$\phantom{\widetilde U=\{\langle y,t\rangle\mid  }\;\&\; F_{t+1}(y)=0\},$

\medskip

$
\widetilde V=\{\langle y,t\rangle\mid (\exists s>t)[ F_s(y)= F_t(y)=1]\;\&\; F_{t+1}(y)=0\}.
$

\medskip

\noindent Suppose now that $\widetilde U\subseteq Z\subseteq \widetilde V$.  We will prove that $F\le_T E\oplus Z$ considering two cases.

\medskip\noindent  {\it Case 1. There are infinitely many elements $\langle y,t\rangle \in \widetilde F\setminus Z.$} Note that for every $y$ there are only finitely many $t$ such that $F_{t}(y)=1$ and $F_{t+1}(y)=0$. Hence, the following $(E\oplus Z)$-c.e. subset of $F$ is infinite:
$$
Y=\{y\mid (\exists t)[\langle y,t\rangle \in \widetilde F\setminus Z]\}\subseteq F.
$$
 Now if $y\in Y$  then for every $s\ge y+1$ and $x<y$ we have
$$  F_{y+1}(x)= F_s(x)=F(x),$$
since otherwise  $  F_s(x)\ne F_{s+1}(x)$ would imply $1= F_t(y)= F(y)= F_s(y)$ and, therefore, $\langle y,t\rangle\in\widetilde U\subseteq Z.$

 Thus, each new $y$ from the infinite  $(E\oplus Z)$-c.e. set $Y$ gives a possibility to  find the value $F(x)$ for each $x<y$. This implies $F\le_T E\oplus Z$.

\medskip \noindent  {\it Case 2. There are only finitely many elements $\langle y,t\rangle \in \widetilde F\setminus Z.$} Then we make only finitely many errors computing $F(x)=\widetilde F(\langle y,0\rangle)$ using the following $(E\oplus Z)$-computable recursive procedure which assumes  $\widetilde F\subseteq Z\subseteq \widetilde V$.

\medskip \noindent  {\it The procedure deciding  whether $\langle y,t\rangle \in \widetilde F$.}
\begin{enumerate}
\item If $\langle y,t\rangle \notin  Z$ then answer ``no''.
\item If $\langle y,t\rangle \in Z$ then $\langle y,t\rangle \in \widetilde V$, and so  $ F_s(y)=1$ for some $s>t$.
\item Due $F=W_e^E=\lim_s  F_s$ there is a $w\ge s$ such that either $y\in W_{e,w}^E$, or $ F_w(y)=1, $  $ F_{w+1}(y)=0.$
\item In the former case answer  ``yes''.
\item In the last case we call recursively the procedure for $\langle y,w\rangle$ by the reduction
$$\langle y,t\rangle \in \widetilde F\iff \langle y,w\rangle \in \widetilde F.$$
Since for each $y$ there are only finitely many $w$ with $ F_w(y)=1, $  $ F_{w+1}(y)=0$ the recursion chain can not be infinite.
\end{enumerate}
Therefore, in both cases  we have proved  $F\le_T E\oplus Z$. 

Consider now the the interval of c.e. sets $U\subseteq V$, where  $U=E\oplus \widetilde U$ and $V=E\oplus\widetilde V.$ Then we have $U\subseteq E\oplus \widetilde F\subseteq V$ for the $(\mathbf e\cup\mathbf f)$-computable set $E\oplus \widetilde F$, and also
$$
U\subseteq X\subseteq V\implies \mathbf e\cup\mathbf f\le \deg(X),
$$
since every set $X$, $U\subseteq X\subseteq V$, must have the form $X=E\oplus Z$, where $\widetilde U\subseteq Z\subseteq \widetilde V$.

\section{The proof of Theorem 5}

 Cooper proved (see~\cite{Co}, 12.3.6 and 12.3.7) that for every $n>1$ there is an $n$-c.e. set $F$ such that $F\not\equiv_T V_e$ for every $(n-1)$-c.e. set $V_e.$ The construction of $F$ can be given via a $\Delta^0_2$-approximation $F(x)=\lim_sF_s(x)$, $F_s(x)\in\{0,1\}$, such that
$$F_0(y)=0\;\&\; {\rm card}\{s\mid F_s(x)\ne F_{s+1}(x) \}\le n$$
for every $x$. Moreover, since  each requirement deals only with one witness at once we have the property
$$x<y<s<t\;\&\; F_s(x)\ne  F_{s+1}(x) \implies  F_s(y)= F_t(y).$$
The last property again holds if we re-define the
 $\Delta^0_2$-approximation for $s=0,1$ by setting $F_0(x)=1$ and $F_1(x)=0$ for each $x$. 
This re-definition also gives
$$
 1\le {\rm card}\{s\mid F_s(x)\ne F_{s+1}(x) \}\le n+1.$$
Then for the sequence of $F$-computable sets
$$
\widetilde F^i=\{\langle y,t\rangle\mid  F_{t}(y)= F(y)\;\&\; F_{t+1}(y)\ne F_t(y)\;\&\; c_t(y)=i\},
$$
where
$$c_t(x)={\rm card}\{s\le t\mid F_s(x)\ne F_{s+1}(x) \}\ge 1,$$
we have

$$\emptyset=\widetilde F^{n+1}=\widetilde F^{n+2}=\widetilde F^{n+3}=\cdots$$
(the equality $\widetilde F^{n+1}=\emptyset$ follows from $F_{t+1}(y)=F(y)$ if $c_t(y)=n+1$, other equalities follow from $c_t(y)\le n+1$), and
$$
F=\{y: \langle y,0\rangle\in \widetilde F^1\}\le_T \widetilde F^1.
$$

Also for every $i$ the set $\widetilde F^i$ is $\widetilde F^{i+1}$-c.e. since
$
\langle y,t\rangle\in \widetilde F^i
$ iff $F_{t+1}(y)\ne F_t(y)$, $c_t(y)=i$, and
$$(\exists s>t)[F_{s+1}(y)\ne F_s(y)\;\&\; c_s(y)=i+1\;\&\; \langle y,s\rangle\notin \widetilde F^{i+1} ].
$$
Similarly with the proof of Theorem 4 we have $U\subseteq \widetilde F\subseteq V$ for
$$
\widetilde F=\bigcup_{i=1}^n\widetilde F^i=\{\langle y,t\rangle\mid F_{t}(y)= F(y)\;\&\; F_{t+1}(y)\ne F_t(y)\}\le_T F,
$$
and the c.e. sets

\medskip

$
U=\{\langle y,t\rangle\mid
(\exists x)(\exists s)[x{<}y{<}s\;\&\;  F_s(x)\ne  F_{s+1}(x)\;\&\;  F_s(y)=F_{t}(y)]
$

$\phantom{ U=\{\langle y,t\rangle\mid  }\;\&\; F_{t+1}(y)\ne F_t(y)\},$

\medskip

$
 V=\{\langle y,t\rangle\mid (\exists s>t)[ F_s(y)= F_t(y)]\;\&\; F_{t+1}(y)\ne F_t(y)\}.
$

\medskip

Suppose now that $ U\subseteq X\subseteq  V$. Due $ \widetilde F^1\equiv_T F$ and $ \widetilde F^{n+1}=\emptyset$  it is enough for  $F\le_T X$ to prove that
$$
\widetilde F^{i+1}\le_T X\implies \widetilde F^{i}\le_T X
$$
for every $i$. Indeed, since $\widetilde F^{i+1}\le_T X$ and $\widetilde F^{i}$ is $\widetilde F^{i+1}$-c.e. the set
$$
Y^i=\{y\mid (\exists t)[\langle y,t\rangle \in \widetilde F^i\setminus X]\}
$$
is $X$-c.e. If $Y^i$ is infinite then we get $\widetilde F^i\le_T F\le_T X$ since each $y\in Y^i$ computes
$$F(x)=F_{y+1}(x)$$
for $x<y$. Indeed,  $\langle y,t\rangle \in \widetilde F^i$
and
$ F_s(x)\ne F_{s+1}(x)$, $s\ge y+1$, would imply $F_t(y)= F(y)= F_s(y)$ and, therefore, $\langle y,t\rangle\in U\subseteq X$.

Let us consider the case  when the set $Y^i$ is finite.  To show $\widetilde F^{i}\le_T X$ in this case we need only to know how to decide whether $\langle y,t\rangle\in \widetilde F^{i}$ for $\langle y,t\rangle\in X$ with $c_t(y)=i$. But if $\langle y,t\rangle\in X$ then $\langle y,t\rangle\in  V$ so that for some $s=w+1>t$ we have
$$
F_{t+1}(y)\ne F_t(y)=F_{w+1}(y)\ne F_{w}(y).
$$
For the least such $s=w+1>t$ we also have $c_w(y)=c_t(y)+1=i+1$ so that
$$
\langle y,t\rangle\in \widetilde F^{i}\iff F_t(y)=F(y)\iff F_w(y)\ne F(y)\iff \langle y,w\rangle\notin \widetilde F^{i+1}.
$$
Hence,  $\widetilde F^{i+1}\le_T X$ implies  $\widetilde F^{i}\le_T X$.

Thus, $U\subseteq X\subseteq  V$ implies $F\le_T X$, and simultaneously $U\subseteq \widetilde F\subseteq  V$ holds for an $F$-computable set $\widetilde F.$

\end{document}